\documentclass[preprint,review,12pt]{elsarticle}

\usepackage{amssymb}

\journal{}

\usepackage{amssymb,latexsym}
\usepackage{amsmath}
\usepackage{graphics}
\usepackage{graphicx}
\newtheorem{theorem}{Theorem}

\newtheorem{definition}{Definition}
\newtheorem{result}{Result}

\newtheorem{remark}{Remark}

\newenvironment{proof}[1][Proof]{\textbf{#1.} }{\ \rule{0.5em}{0.5em}}
\numberwithin{equation}{section} \numberwithin{theorem}{section}
\numberwithin{result}{section} \numberwithin{definition}{section}
\numberwithin{remark}{section}

\begin{document}

\begin{frontmatter}
\title{Admissible Mannheim Curves in Pseudo-Galilean Space $G_3^1$}


\author{M. AKY\.{I}\v{G}\.{I}T, A. Z. AZAK, M. TOSUN}
\address{ Department of Mathematics,
Faculty of Arts and Sciences
\\
Sakarya University, 54187 Sakarya/TURKEY \\}

\begin{abstract}
In this paper, admissible Mannheim partner curves are defined in
pseudo-Galilean space $G_3^1$. Moreover, it is proved that the
distance between the reciprocal points of admissible Mannheim pair
and the torsions of these curves are constant. Furthermore, the
relations between the curvatures and torsions of these curves and
Schell Theorem are obtained. Finally, a result about these curves
being general helix is given.
\end{abstract}

\begin{keyword}
Pseudo-Galilean space, Mannheim curve, admissible curve.

\MSC[2008] 53A35, 51M30
\end{keyword}

\end{frontmatter}

\section{Introduction}\label{S:intro}

In the study of the fundamental theory and the characterizations
of space curves, the corresponding relations between the curves
are the very interesting and important problem. The notion of
Bertrand curves was discovered by J. Bertrand in 1850. He studied
curves in Euclidean 3-space whose principal normals are principal
normals of another curve. In \cite{Liu}, Liu and Wang are
concerned with another kind of associated curves, called Mannheim
curve and Mannheim mate (partner curve).

Liu and Wang have given the definition of Mannheim mate as
follows: Let $\alpha$ and $\alpha ^*$ be two space curves. The
curve $\alpha$ is said to be a Mannheim partner curve of $\alpha
^*$ if there exists a one to one correspondence between their
points such that the principal normal vector of $\alpha$ is
linearly dependent with the binormal vector of $\alpha ^*$ . Also,
they obtained the necessary and sufficient conditions between the
curvature and the torsion for a curve to be the Mannheim partner
curve.

In the elementary differential geometry, there are a lot of works
on Bertrand pair but there are rather a few works on Mannheim
pair. Recently, Bertrand curves and Mannheim curves in
3-dimensional Euclidean and Lorentzian space are studied in
\cite{Ekm} and \cite{Liu,Orb}, respectively.

\section{Preliminaries}\label{S:intro}
The pseudo-Galilean geometry is one of the real Cayley-Klein
geometries (of projective signature (0,0,+,-), explained in
\cite{Mol}). The absolute of the pseudo-Galilean geometry is an
ordered triple $\{ w,f,I\}$ where $w$ is the ideal (absolute)
plane, $f$ is a line in $w$ and $I$ is the fixed hyperbolic
involution of points of $f$. In \cite{Div2, Div3, Mil}, the
pseudo-Galilean space $G_3^1$ has been studied by Divjak.

The group of motions of  $G_3^1$ is a six-parameter group given
(in affine coordinates) by \cite{Mil}

\begin{equation*}
\displaylines{ \overline x  = a + x \cr
 \overline y  = b + cx +
y\cosh \varphi  + z\sinh \varphi  \cr
  \overline z  = d + ex + y\sinh
\varphi  + z\cosh \varphi . \cr}
\end{equation*}

As in \cite{Div3}, pseudo-Galilean scalar product can be written
as
\begin{equation*}
 < v_1 ,v_2  >  = \left\{ \begin{array}{l}
 \,\,\,\,\,\,\,x_1 x_2\,\,\,\, \,\;\;\;\;,\;\;if\;\;x_1  \ne 0\; \vee \;x_2  \ne 0 \\
 y_1 y_2  - z_1 z_2 \;,\;\;if\;\;x_1  = 0\; \wedge \;x_2  = 0 \\
 \end{array} \right.
\end{equation*}
where $v_1  = (x_1 ,y_1 ,z_1 )$ and $v_2  = (x_2 ,y_2 ,z_2 )$. It
leaves invariant the pseudo-Galilean norm of the vector $v =
(x,y,z)$ defined by \cite{Div}
\[
\left\| v \right\| = \left\{ \begin{array}{l}
 \;\;\;\;\;\;\;\;\;\;x\;\;\;\;\;\;\;,\;x \ne 0 \\
 \sqrt {\left| {y^2  - z^2 } \right|} \;\;,\;x = 0. \\
 \end{array} \right.
\]

A vector $v = (x,y,z)$ in $G_3^1$ is said to be non-isotropic if
$x \ne 0$ , otherwise it is isotropic. All unit non-isotropic
vectors are of the form $(1,y,z)$. There are four types of
isotropic vectors: spacelike $(y^2  - z^2  > 0)$, timelike $(y^2 -
z^2  < 0)$ and two types of ligthlike $(y = \pm z)$ vectors. A
non-lightlike isotropic vector is unit vector if $y^2  - z^2  =
\pm 1$  (Figure 1), \cite{Div3}. \\

A trihedron $(T_0 ;e_1 ,e_2 ,e_3 )$, with a proper origin $T_0
(x_0 ,y_0 ,z_0 ) \sim (1:x_0 :y_0 :z_0 )$, is orthonormal in
pseudo-Galilean sense if and only if the vectors $e_1 ,e_2 ,e_3$
have the following form:
\[
e_1  = (1,y_1 ,z_1 ),\;e_2  = (0,y_2 ,z_2 ),\;e_3  =
(0,\varepsilon z_2 ,\varepsilon y_2 )
\]
with $y_2^2  - z_2^2  = \delta$, where each of $\varepsilon
,\delta$ is +1 or -1. Such trihedron $(T_0 ;e_1 ,e_2 ,e_3 )$ is
called positively oriented if $\det (e_1 ,e_2 ,e_3 ) = 1$ holds,
i.e., if $y_2^2  - z_2^2  = \varepsilon$, \cite{Div3}.

Let $\alpha$ be a curve given first by
\begin{equation}
\alpha :I \to G_3^1 ,\,\,\alpha (t) = (x(t),y(t),z(t))
\end{equation}
where $I \subseteq \mathbb{R}$ and $x(t),y(t),z(t) \in C^3$ (the
set of three-times continuously differentiable functions) and  $t$
run through a real interval. A curve $\alpha$   given by (2.1) is
admissible if $x'(t)$ is not equal to zero, \cite{Div3}.

The curves in pseudo-Galilean space   are characterized as follows
\cite{Div2,Erj}.

An admissible curve $\alpha$ in $G_3^1$ can be parameterized by
arc length $t = s$, given in coordinate form

\begin{equation}
\alpha (s) = (s,y(s),z(s)).
\end{equation}

For an admissible curve $\alpha :I \to G_3^1 ,\,I \subseteq
\mathbb{R}$, the curvature $\kappa _\alpha  (s)$ and the torsion
$\tau _\alpha (s)$ are defined by

\begin{equation}
\begin{array}{l}
 \kappa _\alpha  (s) = \sqrt {\left| {y''^2 (s) - z''^2 (s)} \right|}  \\
 \tau _\alpha  (s) = \frac{{\det (\alpha '(s),\alpha ''(s),\alpha '''(s))}}{{\kappa _\alpha  ^2(s) }}. \\
 \end{array}
\end{equation}

The associated trihedron is given by

\begin{equation}
\begin{array}{l}
 T_\alpha  (s) = \alpha '(s) = (1,y'(s),z'(s)) \\
 N_\alpha  (s) = \frac{1}{{\kappa _\alpha  (s)}}\alpha ''(s) = \frac{1}{{\kappa _\alpha  (s)}}(0,y''(s),z''(s)) \\
 B_\alpha  (s) = \frac{1}{{\kappa _\alpha  (s)}}(0,z''(s),y''(s)). \\
 \end{array}
\end{equation}

The vectors $T_\alpha  ,N_\alpha  ,B_\alpha$ are called the
vectors of tangent, principal normal and binormal line of
$\alpha$, respectively. The curve $\alpha$ given by (2.2) is
timelike if $N_\alpha  (s)$ is a spacelike vector. For derivatives
of the tangent (vector) $T_\alpha$, principal normal $N_\alpha$
and binormal
 $B_\alpha$, respectively, the following Frenet formulas hold

\begin{equation}
\begin{array}{l}
 T_\alpha  '(s) = \kappa _\alpha  (s)N_\alpha  (s) \\
 N_\alpha  '(s) = \tau _\alpha  (s)B_\alpha  (s) \\
 B_\alpha  '(s) = \tau _\alpha  (s)N_\alpha  (s). \\
 \end{array}
\end{equation}

\section{Admissible Mannheim Curves in Pseudo-Galilean Space $
G_3^1$}\label{S:intro}

In \cite{Liu}, Liu and Wang have studied Mannheim partner curves
in 3-dimensional Euclidean space $E^3$. In this section, we
introduce the notion of Mannheim partner curves in $G_3^1$ for
admissible curves in the following way:

\begin{definition}

 Let $\alpha$ and $\alpha ^*$ be two admissible curves
with non-zero $\kappa _\alpha  (s)$, $\kappa _\alpha ^* (s)$,
$\tau _\alpha (s)$, $\tau _\alpha ^* (s)$ for each $s \in I$ and
$\{ T_\alpha ,N_\alpha ,B_\alpha \}$ and ${\kern 1pt} \{ T_\alpha
^* ,N_\alpha ^* ,B_\alpha ^* \}$ be Frenet frame in $G_3^1$ along
$\alpha$ and $\alpha ^*$, respectively. If there exists a
corresponding relationship between the admissible curves $\alpha$
and $\alpha ^*$, such that at the corresponding points of the
admissible curves, principal normal lines $N_\alpha$ of $\alpha$
coincides with the binormal lines $B_\alpha ^*$ of $\alpha ^*$,
then $\alpha$ is called an admissible Mannheim curve and $\alpha
^*$ an admissible Mannheim partner curve of $\alpha$. The pair
$(\alpha ,\alpha ^* )$ is said to be an admissible Mannheim pair
in $G_3^1$ (see Figure 2).
\end{definition}

Throughout this paper we assume that $\alpha$ is timelike curve
and its Mannheim partner $\alpha ^*$ is spacelike curve, i.e.,
$N_\alpha$, the principal normal of $\alpha$ curve, is spacelike
vector and $N_\alpha ^*$, the principal normal of $\alpha ^*$
curves, is timelike vector.

\begin{theorem}\label{T:2.1}Let $(\alpha ,\alpha ^* )$ be an admissible Mannheim
pair in $G_3^1$. Then function $\lambda$ is constant.
\end{theorem}
\begin{proof}
Let $(\alpha ,\alpha ^* )$ be an admissible Mannheim pair and the
Frenet frame be $\{ T_\alpha  ,N_\alpha  ,B_\alpha  \}$ and
${\kern 1pt} \{ T_\alpha ^* ,N_\alpha ^* ,B_\alpha ^* \}$ on
$G_3^1$ along $\alpha$ and $\alpha ^*$, respectively. Let the
point of $\alpha ^*$ be corresponding to the point of $\alpha$,
then the admissible Mannheim curve $\alpha ^*$ is given by

\begin{equation}
\alpha ^* (s) = \alpha (s) + \lambda (s)N_\alpha  (s).
\end{equation}

\noindent By taking the derivative of equation (3.1) with respect
to $s$ and applying the Frenet formulas, we get

\begin{equation}
T_\alpha ^* (s)\frac{{ds^* }}{{ds}} = T_\alpha  (s) + \lambda
'(s)N_\alpha  (s) + \lambda (s)\tau _\alpha  (s)B_\alpha  (s).
\end{equation}

\noindent Since $N_\alpha$ is coincident with $B_\alpha ^*$ in the
same direction, we have

\begin{equation}
\lambda '(s) = 0.
\end{equation}

\noindent This gives $\lambda (s) =$ constant.
\end{proof}

In other words, the theorem proves that the distance between
corresponding points of curve $\alpha$ and its Mannheim pair
$\alpha ^*$ is constant.

\noindent It should be noted that the relationship between a curve
and its Mannheim pair is a reciprocal one, i.e., if a curve
$\alpha ^*$ is a Mannheim pair of a curve $\alpha$, then  $\alpha$
is also Mannheim pair of $\alpha ^*$.

\begin{theorem}\label{T:2.1}
Let $\alpha$  be an admissible curve with arc length parameter
$s$. If $\alpha$ is an admissible Mannheim curve, then the torsion
$\tau _\alpha$ of admissible curve $\alpha$ is constant. (Converse
of the theorem is also true.)
\end{theorem}

\noindent \begin{proof} Let $(\alpha ,\alpha ^* )$ be an
admissible Mannheim curves. Then we have

\begin{equation}
\begin{array}{l}
 T_\alpha  (s) = \cosh \theta T_\alpha ^* (s) - \sinh \theta N_\alpha ^* (s) \\
 B_\alpha  (s) =  - \sinh \theta T_\alpha ^* (s) + \cosh \theta N_\alpha ^* (s) \\
 \end{array}
\end{equation}
where $\theta$ is the angle between $T_\alpha$ and $T_\alpha ^*$
at the corresponding points of $\alpha$ and $\alpha ^*$ (see
Figure 2). On the other hand, from the last equation, we get
\begin{equation}
\begin{array}{l}
 T_\alpha ^* (s)\,\,\, = \cosh \theta T_\alpha  (s) + \sinh \theta B_\alpha  (s) \\
 N_\alpha ^* (s) = \sinh \theta T_\alpha  (s) + \cosh \theta B_\alpha  (s). \\
 \end{array}
\end{equation}

\noindent By taking the derivative of the equation (3.5) with
respect to $s$ we get

\begin{equation}
\begin{array}{l}
 \tau _\alpha ^* (s)B_\alpha ^* (s)\frac{{ds^* }}{{ds}} = \frac{{d(\sinh \theta )}}{{ds}}T_\alpha  (s) + \sinh \theta \kappa _\alpha  (s)N_\alpha  (s) \\
 \,\,\,\,\,\,\,\,\,\,\,\,\,\,\,\,\,\,\,\,\,\,\,\,\,\,\,\,\,\,\,\,\,\,\,\,\,\,\,\,\,\, + \cosh \theta \tau _\alpha  (s)N_\alpha  (s) + \frac{{d(\cosh \theta )}}{{ds}}B_\alpha  (s). \\
 \end{array}
\end{equation}

\noindent Since the principal normal vector $N_\alpha$ of the
curve $\alpha$ and the binormal vector $B_\alpha ^*$ of its
Mannheim partner curve is linearly dependent, from the last
expression, we reach

\begin{equation*}
\theta  =constant.
\end{equation*}

\noindent By considering equation (3.2), it follows that

\begin{equation}
T_\alpha ^* (s)\frac{{ds^* }}{{ds}} = T_\alpha  (s) + \lambda \tau
_\alpha  (s)B_\alpha  (s).
\end{equation}

\noindent If we take into consideration equations (3.5) and (3.7),
then we obtain

\begin{equation}
\lambda \tau _\alpha  (s) \coth \theta  = 1.
\end{equation}

\noindent Taking $u = \lambda \coth \theta$ and using Theorem 3.1,
we reach

\begin{equation}
\tau _\alpha   (s) = \frac{1}{u}.
\end{equation}

\noindent The fact that $\tau_\alpha$ is a constant and this
completes the proof.
\end{proof}

\begin{theorem}
(Schell's Theorem). Let $(\alpha ,\alpha ^* )$ be an admissible
Mannheim pair in $G_3^1$. The product of torsions $\tau _\alpha$
and $\tau _\alpha ^*$ at the corresponding points of the
admissible Mannheim curves is constant, where the torsion $\tau
_\alpha$ belong to $\alpha$ and the torsion $\tau _\alpha ^*$
belong to $\alpha ^*$.
\end{theorem}

\noindent \begin{proof} If we take $\alpha ^*$ instead of
$\alpha$, then we can rewrite the equation (3.1) as follow:

\begin{equation}
\alpha (s) = \alpha ^* (s) - \lambda B_\alpha ^* (s)
\end{equation}

\noindent By taking derivative of equation (3.10) with respect to
$s$ and using equation (3.4), we get

\begin{equation}
\tau _\alpha ^* (s) =  \frac{1}{\lambda }\tanh \theta .
\end{equation}

\noindent Multiplying both sides of equation (3.8) by the
corresponding sides of equation (3.11), we have

\begin{equation}
\tau _\alpha  (s)\tau _\alpha ^* (s) = \frac{{\tanh ^2 \theta
}}{{\lambda ^2 }} =constant.
\end{equation}

\noindent Hence the proof is completed.
\end{proof}

\begin{theorem}
Let $(\alpha ,\alpha ^* )$ be an admissible Mannheim pair in
$G_3^1$ and $\kappa _\alpha  ,\tau _\alpha,\kappa _\alpha ^* ,\tau
_\alpha ^*$ be a curvatures and torsions of $\alpha$ and $\alpha
^*$, respectively. Then their curvatures and torsions satisfy the
following relations:

\begin{equation*}
\begin{array}{l}
 i)\;\;\;\kappa _\alpha ^* (s) = -\frac{{d\theta }}{{ds^*}} \\
 ii)\;\;\kappa _\alpha  (s) =  \tau _\alpha ^* (s)\sinh \theta \frac{{ds^* }}{{ds}} \\
 iii)\;\;\tau _\alpha  (s) = -\tau _\alpha ^* (s)\cosh \theta \frac{{ds^* }}{{ds}}. \\
 \end{array}
\end{equation*}
\end{theorem}

\noindent  \begin{proof} \textit{i)} If we consider equation
(3.4), we get
\begin{equation}
 < T_\alpha  (s),T_\alpha ^* (s) >  = \cosh \theta .
\end{equation}
By taking the derivative of the last equation with respect to
$s^*$ and using the Frenet formula for $\alpha$ and $\alpha ^*$,
we reach
\begin{equation}
 < \kappa _\alpha  (s)N_\alpha  (s)\frac{{ds}}{{ds^* }},T_\alpha ^* (s) >  +  < T_\alpha  (s),\kappa _\alpha ^* (s)N_{_\alpha  }^* (s) >  = \sinh \theta \frac{{d\theta }}{{ds^* }}
\end{equation}
By considering the principal normal vector $N_\alpha$ of $\alpha$
and binormal vector $B_\alpha ^*$ of $\alpha ^*$ are linearly
dependent and using equations (3.4) and (3.14), we obtain
\begin{equation}
\kappa _\alpha ^* (s) = -\frac{{d\theta }}{{ds^*}}.
\end{equation}

If we consider the equations (2.5), (3.4), (3.5) and the scalar
products of $<T_\alpha ,B_\alpha ^*>$, $ < B_\alpha ,B_\alpha ^*
>$, we can easily prove the items $ii),\,\,iii)$
 of the theorem, respectively.
\end{proof}

From equations $ii)$ and $iii)$ of the last theorem, we get
\begin{equation*}
\frac{{\kappa _\alpha  (s)}}{{\tau _\alpha  (s)}} =  - \tanh
\theta  =constant.
\end{equation*}
Hence, we give following result.
\begin{result} If $\theta$ is
constant where $\theta$ is angle between $T_\alpha$ and $T_\alpha
^*$   at the corresponding points of $\alpha$ and $\alpha ^*$. An
admissible Mannheim curve $\alpha$ in $G_3^1$ is a general helix.

\end{result}

We give the following result (helices as Mannheim partner curves).

\begin{result} Let $(\alpha ,\alpha ^* )$ be an admissible
Mannheim pair with arc length parameter in $G_3^1$. If $\alpha$ is
a generalized helix, then $\alpha ^*$ is a straight line,
\cite{Liu}.
\end{result}

\begin{remark}
The method used above can be considered for spacelike curve
$\alpha$ and its Mannheim partner timelike curve $\alpha ^*$. So
we can easily obtain the similar theorems and results for these
curves as above theorems and results.
\end{remark}

\newpage
Figure 1. The points in pseudo-Galilean space

Figure 2. The admissible Mannheim partner curves

\end {document}
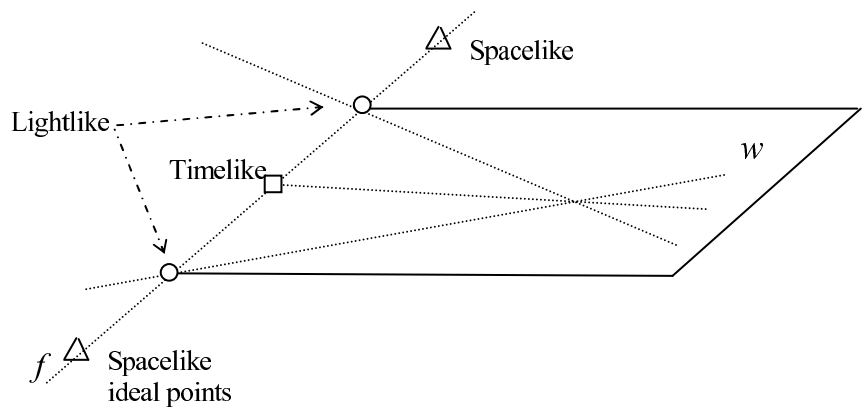
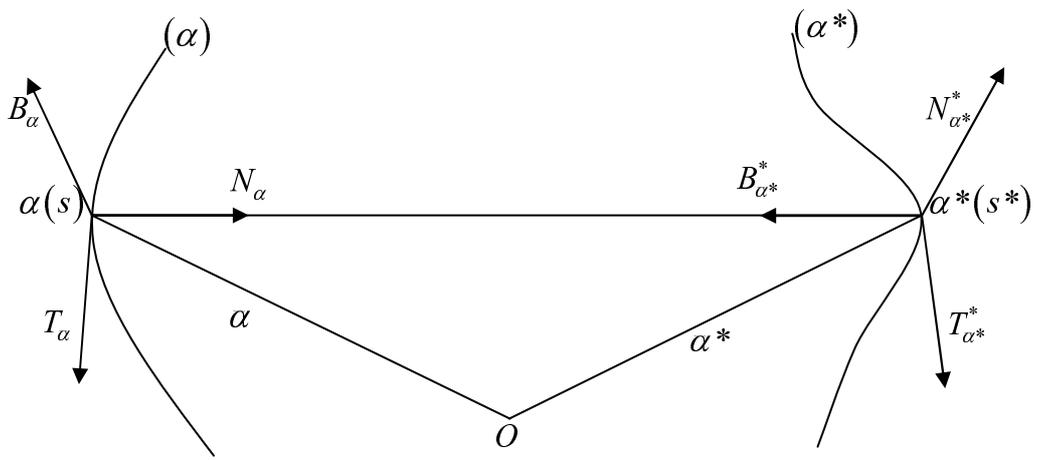
